\documentclass[12pt]{amsart}
\usepackage{geometry,a4wide,amssymb,amsfonts,amsrefs,amsmath}                
\usepackage{bbm,xcolor}
\usepackage{dsfont}
\usepackage{yfonts}
\usepackage[T1]{fontenc}
\usepackage{mathtools}
\usepackage{mathrsfs}
\usepackage{hyperref}  
\usepackage{footmisc} 

\newcommand{\Z}{\mathbb{Z}}
\newcommand{\h}{\mathtt{h}}
\newcommand{\A}{\mathcal{A}}
\newcommand{\N}{\mathbb{N}}
\newcommand{\bL}{\mathbb{L}}
\newcommand{\R}{\mathbb{R}}
\newcommand{\F}{\mathcal{F}}
\newcommand{\s}{\mathfrak{S}}
\newcommand{\0}{\mathtt{0}}
\newcommand{\M}{\mathcal{M}}
\newcommand{\G}{\mathscr{G}}
\newcommand{\cov}{\mathtt{cov}}

\newcommand{\La}{\mathcal{L}}
\newcommand{\1}{\mathds{1}}

\theoremstyle{plain}

\newtheorem{thm}{Theorem}[section]
\newtheorem{lem}[thm]{Lemma}
\newtheorem{prop}[thm]{Proposition}

\newtheorem{defin}[thm]{Definition}

\newtheorem{remark}[thm]{Remark}

  \newtheorem{manualtheoreminner}{Theorem}
\newenvironment{manualtheorem}[1]{%
  \renewcommand\themanualtheoreminner{#1}%
  \manualtheoreminner
}{\endmanualtheoreminner}
  
\title[Eigenfunctions for the Dyson model in a field]{The eigenfunctions of the transfer operator for the Dyson model in a field}
\author{Mirmukhsin Makhmudov}
\address{Mathematical Institute,  Leiden University, Einsteinweg 55, 2333 CC Leiden, The Netherlands}
\email{m.makhmudov@math.leidenuniv.nl}

\begin{document}
\begin{abstract}
The recent works \cite{EFMV2024} and \cite{JOP2023} have studied the spectral properties of the Dyson model in the absence of an external field.
This paper is a continuation of  \cite{EFMV2024} and aims to bridge the gap in the literature by investigating the Dyson model in a field.

\noindent
In this paper, we prove that, for high temperatures or strong magnetic fields, there exists
a non-negative, integrable (with respect to the unique half-line Gibbs measure) eigenfunction of the transfer operator for the Dyson model if $\alpha\in(\frac 3 2,2]$. 
However, unlike in the zero-magnetic-field case, this eigenfunction is not continuous.
\end{abstract}
\maketitle

\section{Introduction}

The study of the so-called \textit{equilibrium states}, specific types of invariant measures, has been a central topic in the ergodic theory of dynamical systems since the 1970s. The foundations of this theory were laid in the seminal works of Bowen and Ruelle at the end of the decade.
Equilibrium states also play a crucial role in equilibrium statistical mechanics. A pivotal result in this context is the classical variational principle, proved by Lanford and Ruelle in 1969, which establishes that, under certain conditions, translation-invariant Gibbs measures coincide with equilibrium states.
Recall that an equilibrium state $\mu$ for a continuous potential $\phi$ on the one-sided shift space (we also refer to it as the half-line shift space) $X_+:=E^{\Z_+}$, $|E|<\infty$ are probability measures on $X_+$ which are invariant under the shift (translation) map $S: X_+\to X_+$, $(Sx)_j=x_{j+1}$, $j\in\Z_+$ and solves the variational equation:
\begin{equation}
    \int_{X_+}\phi d\mu +\mathfrak{h}(\mu)
    =
    \sup
    \Big\{
    \int_{X_+}\phi d\tau +\mathfrak{h}(\tau): \; \tau\in\M_{1,S}(X_+)
    \Big\},
\end{equation}
where $\M_{1,S}(X_+)$ is the simplex of translation-invariant probability measures on $X_+$ and $\mathfrak{h}(\tau)$
denotes the measure-theoretic entropy of the translation-invariant measure $\tau$.

A natural and effective approach to studying equilibrium states for a potential $\phi$
is through the transfer operator $\La_\phi$, which acts on the space of functions on $X_+$, particularly on continuous functions, via:
\begin{equation}
    \La_\phi f(x):=\sum_{a\in E}e^{\phi(ax)} f(ax),\;\; f\in \R^{X_+},\;\; x\in X_+.
\end{equation}
In fact, if $\h$ is an eigenfunction of $\La_\phi$, i.e., $\La_\phi \h = \lambda \h$, and $\nu$ is an eigenprobability for its adjoint $\La^*_\phi: C(X_+)^*\to C(X_+)^*$, i.e., $\La^*_\phi \nu = \lambda \nu$, where 
$\lambda$ is the spectral radius of $\La_\phi$, then the measure $d\mu =\h\cdot d\nu$ is an equilibrium state for $\phi$ and vice versa (see Proposition 3.1 in \cite{EFMV2024}).
Classical fixed-point theorems ensure the existence of an eigenprobability for the adjoint $\La_\phi^*$, however, the existence of an eigenfunction for $\La_\phi$ is not always guaranteed. 
This issue has been extensively studied by D. Ruelle, P. Walters and others under various regularity conditions on $\phi$ \cites{Ruelle-book, Walters1975, Walters1978, Walters2001}. 
Notably, every potential $\phi$ satisfying these regularity conditions shares two key properties:
\begin{itemize}
    \item[(P1)] for all $\beta > 0$, the potential $\beta \phi$ admits a unique equilibrium state;
    \item[(P2)] the transfer operator $\La_\phi$ is stable under "smooth" perturbations of $\phi$ in the sense that if $\La_\phi$ has a (principal) eigenfunction, then $\La_{\phi+\upsilon}$ also has a (principal) eigenfunction with the same regularity properties for every local function $\upsilon: X_+\to\R$.
\end{itemize}
Note that the first property rules out the possibility of phase transitions, i.e., the existence of multiple equilibrium states for some $\beta$, a phenomenon of significant interest in statistical mechanics. Consequently, these classical conditions on $\phi$ exclude important one-dimensional models.
Efforts to extend the theory to more general one-dimensional models have emerged in recent studies. 
In 2019, Cioletti et al. analyzed transfer operators for product-type potentials, which fall outside existing classes but still satisfy (P1) and (P2) \cite{CDLS2017}.
A more notable development came in 2023 when Johansson, Öberg, and Pollicott analyzed the so-called Dyson potential, the potential corresponding to the 1D long-range Ising model \cite{JOP2023} which is defined by:
\begin{equation}\label{eq: Dyson potential on half-line}
    \phi( x) = hx_0 + \beta\sum_{n=1}^\infty \frac {x_0x_n}{n^\alpha},\quad x\in X_+:=\{\pm 1\}^{\Z_+}.
\end{equation}
where $\alpha>1$ is the decay rate of the coupling, $\beta \geq 0$ is the inverse temperature, and $h \in \mathbb{R}$ represents the strength of the external field.
Shortly thereafter, the author of this paper, in collaboration with van Enter, Fernández, and Verbitskiy, introduced a new approach to addressing the eigenfunction problem for long-range models \cite{EFMV2024}.
Following the developments in \cites{EFMV2024, JOP2023}, the following result on the Dyson model is known:

\begin{thm}\cites{EFMV2024, JOP2023}\label{our and JOP's result when h=0}
Let $\phi$ be the Dyson potential (\ref{eq: Dyson potential on half-line}). 
Suppose $h=0$,
then:
\begin{itemize}
    \item[(i)] for all $\alpha>1$ and for $\beta\geq 0$ sufficiently small, the unique equilibrium state $\mu_+$ for $\phi$ is equivalent to the unique half-line Gibbs state $\nu$ for $\phi$, i.e., $\mu_+\ll\nu$ and $\nu\ll\mu_+$.
    In particular, the Perron-Frobenius transfer operator $\La_{\phi}$ has an eigenfunction in $L^1(X_+,\nu)$;  
    \item[(ii)] furthermore, if $\alpha>\frac{3}{2}$, then for all $0<\beta<\beta_c$
    there exists a continuous version
    of the Radon-Nikodym density $\frac{d\mu_+}{d\nu}$, and thus the Perron-Frobenius transfer operator $\La_{\phi}$ has a positive continuous eigenfunction, where $\beta_c$ is the critical temperature of the phase transitions for the standard Dyson model (see Theorem \ref{phase diagram of Dyson model} in the next section). 
\end{itemize}
\end{thm}
However, both studies, \cite{JOP2023} and 
\cite{EFMV2024}, do not cover the case of the Dyson potential with a nonzero external field $h\neq 0$.
The approach in \cite{JOP2023} relies heavily on the random cluster representation of the Dyson model. 
The central obstacle to extending the method in \cite{JOP2023} to non-zero external fields is the loss of symmetry, essential for the random cluster representation, which disrupts the cluster decay analysis.
The method developed in \cite{EFMV2024} requires a certain sum of two-point functions to be uniformly bounded, a condition that fails for the Dyson model in a field.
In this paper, we adapt the technique from \cite{EFMV2024} to treat the Dyson model in a nonzero field. Our main result is as follows:

\begin{manualtheorem}{A}\label{main result}

Suppose $\alpha\in(\frac{3}{2}, 2]$, $\beta\geq 0$ and 
$|h|>0$ is sufficiently large $(|h|>2\beta\zeta(\alpha) + \log 4\beta\zeta(\alpha)
\text{ is enough, here } \zeta \text{ is the Riemann zeta function})$.
Then 
\begin{itemize}
    \item[(i)] the Dyson potential $\phi$ has a unique equilibrium state $\mu_+$ and there also exists a unique eigenprobability $\nu$ of $\La_\phi^*$;   
    \item[(ii)] $\mu_+$ is absolutely continuous with respect to 
    $\nu$, i.e., $\mu_+\ll \nu$.
    In particular, the Perron-Frobenius transfer operator $\La_\phi$ admits an integrable eigenfunction corresponding to the spectral radius $\lambda=e^{P(\phi)}$.
    \item[(iii)] The Radon-Nikdoym derivative $\frac{d\mu_+}{d\nu}$ does not have a continuous version.
    In particular, the Perron-Frobenius transfer operator does not have a continuous principal eigenfunction.
\end{itemize}
\end{manualtheorem}

\begin{remark}
\begin{itemize}
    \item[(i)] Theorem \ref{main result} is proven under a strong uniqueness condition for the Dyson interaction $\Phi$ (see (\ref{eq: Dyson interaction}) and Section \ref{Spec section}), a concept introduced in \cite{EFMV2024}. 
    This condition is implied by the Dobrushin uniqueness condition (DUC) on $\Phi$ (Section \ref{DUC section}). 
    The interaction $\Phi$ satisfies DUC if either 
    1) $\beta$ is sufficiently small so that $\beta\leq \frac{1}{2\zeta(\alpha)}$ or 
    2) $|h|$ is sufficiently large so that $|h|>2\beta\zeta(\alpha) + \log 4\beta\zeta(\alpha)$.
    Note that in the first version of the Dobrushin uniqueness condition, there is no constraint on the external field $h$.
    We prove Theorem \ref{main result} for large $|h|$'s, however, the statement of Theorem \ref{main result} remains valid and the same proof works for the first version of the DUC which covers the domain $0\leq \beta <\frac{1}{2\zeta(\alpha)}$ and $h\neq 0$.

\noindent
    We, in fact, believe that Theorem \ref{main result} is true for all $h\in\R$ but $h=0$, i.e., the transfer operator $\La_\phi$ has an integrable but not continuous eigenfunction for all $\beta\geq 0$, $h\neq 0$.
    \item[(ii)] 
    It is worth noting that Theorem \ref{main result} suggests a certain similarity between the Dyson potential in a field and product-type potentials \cite{CDLS2017}.
    Note that for a product-type potential, there exists a unique equilibrium state $\mu_+$ and eigenprobability $\nu$ which both are Bernoulli measures such that $\nu=\prod_{n=0}^\infty \lambda_n$ with $\lambda_n(1)=p_n=\frac{\exp{(\beta\sum_{i=1}^ni^{-\alpha})}}{2\cosh{(\beta\sum_{i=1}^ni^{-\alpha})}}$ and $\mu_+=\prod_{n=0}^\infty \lambda_\infty$ \cites{CDLS2017, EFMV2024}. 
    Thus 
    one can check that
    $$
    \int_E\sqrt{\frac{d\lambda_\infty}{d\lambda_n}} d\lambda_n
    =
    \sqrt{p_\infty p_n}+\sqrt{(1-p_\infty)(1-p_n)}
    =1+O^*((p_n-p_\infty)^2),
    $$
    therefore, 
    \begin{equation}\label{input eq to kakutani th}
         \log \int_E\sqrt{\frac{d\lambda_\infty}{d\lambda_n}} d\lambda_n
    =
    O^*((p_n-p_\infty)^2)=O^*(n^{-2\alpha+2}),
    \end{equation}
    here $u(t)=O^*(v(t))$ means the existence of $C_1, C_2>0$ such that $C_1 |u(t)|\leq |v(t)|\leq C_2|u(t)|$ for all sufficiently small $t$.
    Then by (\ref{input eq to kakutani th}), Kakutani's theorem implies that $\mu_+\ll \nu$ if $\alpha>  \frac{3}{2}$ and $\mu_+\perp \nu$ if $\alpha\leq \frac{3}{2}$.
    
    In the case of the Dyson potential, in the presence of a strong external field, spin-spin correlations become negligible, as they decay rapidly.
    Consequently, individual spins behave like independent random variables.
    Thus, we expect a phenomenon similar to that observed with product-type potentials, as described above, to occur in the Dyson model under an external field: for $\alpha \in \left(\frac{3}{2}, 2\right]$, we have $\mu_+ \ll \nu$, while for $\alpha \in \left(1, \frac{3}{2}\right]$, $\mu_+ \perp \nu$.
\end{itemize}
\end{remark}
Note that the potential $\phi_{\mathtt{0}}(x)=\beta\sum_{n=1}^\infty\frac{x_0 x_n}{n^\alpha},\; x\in X_+=\{\pm 1\}^{\Z_+}$
violates both (P1) and (P2) properties mentioned above.
In fact, for all $\alpha\in(1,2]$, there exists $\beta_{c}(\alpha)>0$ such that for all $\beta>\beta_c(\alpha)$, $\phi_{\0}$ has multiple equilibrium states \cites{Dyson1969, FS1982}. 
As regards the second property, 
the second parts of Theorem \ref{our and JOP's result when h=0} and Theorem \ref{main result} yields that the transfer operator $\La_{\phi_{\mathtt{0}}}$ 
is unstable under "smooth" perturbations of $\phi_{\mathtt{0}}$ since for $\alpha\in (\frac{3}{2}, 2]$ and $\beta<\beta_c$, the principal eigenfunction of $\La_{\phi_{\mathtt{0}}}$ is continuous 
whereas for $h\neq 0$, $\La_{\phi_{\mathtt{0}}+h\sigma_0}$ admits an integrable but not continuous eigenfunction, where $\sigma_0(x)=x_0$ for $x\in X_+$.

The paper is organized as follows.
In the Section \ref{Spec section}
we define the notions of specifications and the interactions. 
Here we also introduce the interaction of the Dyson model.
Section \ref{DUC section} is dedicated to the Dobrushin uniqueness condition and its consequences which are important for this paper.
In Section \ref{Section: Construction of intermediate interactions}, we discuss the construction of the intermediate interactions which are instrumental in the proof of Theorem \ref{main result}.
Section \ref{Section: Proof of main result} is dedicated to the proof of Theorem \ref{main result}.

\section{Specifications and Interactions}\label{Spec section}

In this paper, we always consider the ferromagnetic spin space $E=\{\pm 1\}$.
For a countable set $\bL$ of sites, we consider a configuration space $\Omega=E^\bL$.
For $\bL$, we either consider the set of integers $\Z$, or the set of negative integers $-\N$, or the set of non-negative integers $\Z_+$.
As a product of compact sets, $\Omega$ is compact and also metrizable topological space.
In fact, for any $\theta\in (0,1)$, $d(x,y)=\theta^{\mathbf{n}(x,y)}$ metrizes $X_+$, here for $x, y\in X_+$, $\mathbf{n}(x,y)=\min\{j\in\Z_+: x_j\neq y_j\}$.
$\F$ denotes the Borel sigma-algebra of $X$.
For any $\Lambda\subset \bL$, a sub sigma-algebra $\F_\Lambda\subset\F$ is defined as the minimal sigma-algebra so that for all $j\in\Lambda$, the functions $\sigma_j: \Omega\to E$ are measurable where  $\xi\in \Omega\mapsto\sigma_j(\xi)=\xi_j\in\{\pm1\}$.
For configurations $\xi, \eta\in\Omega$ and a volume $\Lambda\subset \bL$, $\xi_\Lambda\eta_{\Lambda^c}$ denotes the concatenated configuration, i.e., $(\xi_\Lambda\eta_{\Lambda^c})_i=\xi_i$ if $i\in\Lambda$ and $(\xi_\Lambda\eta_{\Lambda^c})_i=\eta_i$ if $i\in\Lambda^c=\bL\setminus\Lambda$.

A \textit{specification} $\gamma$, a regular family of conditional probabilities, is defined as a \textit{consistent} family of proper probability kernels $\gamma_\Lambda:\F\times \Omega\to [0,1]$ 
indexed by finite subsets $\Lambda$ of $\bL$, denoted by $\Lambda\Subset\bL$ \cite{Georgii-book}*{Chapter 1}.
A probability measure $\mu\in\M_{1}(\Omega)$ is \textit{Gibbs} for a specification $\gamma$ if it is consistent with $\gamma$, i.e., for all $\Lambda\Subset\bL$, $\mu=\mu\gamma_\Lambda$, here
$B\in\F\mapsto \mu\gamma_\Lambda(B):=\int_\Omega\gamma_\Lambda(B|\eta)\mu(d\eta)$.
The set of Gibbs measures for the specification $\gamma$ is denoted by $\G(\gamma)$.

The classical approach to specifications is via interactions.
An interaction $\Phi$ is a family $(\Phi_\Lambda)_{\Lambda\Subset\bL}$ of functions $\Phi_\Lambda: \Omega\to \R$ such that $\Phi_\Lambda\in\F_\Lambda$, i.e., each $\Phi_\Lambda$ is independent of coordinates outside $\Lambda$.
An interaction is called \textit{Uniformly Absolutely Convergent} (UAC) if for all $j\in \bL$,
$$
\sum_{j\in V\Subset\bL} \sup_{\omega\in\Omega}|\Phi_V(\omega)|<\infty.
$$
For a UAC interaction $\Phi$, the Hamiltonian in a volume $\Lambda\Subset\bL$ is 
$H_\Lambda:=\sum_{\substack{V\cap\Lambda\neq \emptyset}}\Phi_V$ and
the associated specification (density) $\gamma^\Phi$ is defined by
\begin{equation}\label{Boltzmann ansatz}
    \gamma^\Phi_\Lambda(\omega_\Lambda|\omega_{\Lambda^c})
    :=
    \frac{e^{-H_\Lambda(\omega)}}{Z_{\Lambda}(\omega_{\Lambda^c})},
    \;\;
    \omega\in\Omega,
\end{equation}
where $Z_\Lambda$ is the normalization constant, also known as the partition function, which is $Z_\Lambda(\omega_{\Lambda^c}):=\sum_{\xi_{\Lambda}\in E^{\Lambda}}e^{-H_\Lambda(\xi_\Lambda\omega_{\Lambda^c})}$.
Note that the specification $\gamma^\Phi=(\gamma_\Lambda^\Phi)_{\Lambda\Subset\bL}$ is restored from (\ref{Boltzmann ansatz}) by 
$$
\gamma^\Phi_\Lambda(B|\omega):=\sum_{\xi_\Lambda\in E^\Lambda}\1_B(\xi_\Lambda\omega_{\Lambda^c})\gamma_\Lambda^\Phi(\xi_\Lambda|\omega_{\Lambda^c}), \;\; B\in\F,\; \omega\in\Omega.
$$
This paper focuses on the Dyson model, a cornerstone of one-dimensional statistical mechanics renowned for its long-range interactions and critical phenomena.
It is defined on the lattice $\bL=\Z$ by 
\begin{equation}\label{eq: Dyson interaction}
    \Phi_{\Lambda}(\omega):=\begin{cases} 
      -\frac{ \beta\omega_i\omega_j}{|i-j|^{\alpha}}, & \text{ if }\Lambda=\{i,j\}\subset\Z, \ i\neq j; \\
      -h \omega_i, & \Lambda=\{i\};\\
      0, & \text{otherwise},
   \end{cases}
\end{equation}
where $\alpha\in(1,2]$ is a parameter describing the decay of the interaction strength, $\beta\geq 0$ is the inverse temperature and $h\in\R$ represents the external field.
In fact, the Dyson potential defined in (\ref{eq: Dyson potential on half-line}) is related to the Dyson interaction $\Phi$ by
$$
\phi=-\sum_{0\in\Lambda\Subset\Z_+}\Phi_\Lambda.
$$
The following theorem about the phase diagram of the Dyson model is well-known.
\begin{thm}\label{phase diagram of Dyson model}\cites{Dyson1969,FS1982, FV-book} Suppose $\Phi$ be the Dyson interaction given in (\ref{eq: Dyson interaction}).
    \begin{itemize}
        \item[(i)] If $h\neq 0$, then $\Phi$ has a unique Gibbs measure, i.e., $|\G(\gamma^\Phi)|=1$;
        \item[(ii)] If $h=0$, there exists a finite critical temperature $\beta_c(\alpha)>0$ such that for all $\beta\in [0,\beta_c(\alpha))$, $|\G(\gamma^\Phi)|=1$ and for $\beta>\beta_c(\alpha)$, $\Phi$ exhibits phase transitions, i.e., $|\G(\gamma^\Phi)|>1$.
    \end{itemize}
\end{thm}

\section{The Dobrushin uniqueness condition and its corollaries}\label{DUC section}
Dobrushin uniqueness condition is one of the most general criteria for the uniqueness of Gibbs states.
We discuss it in the framework of a general countable set $\bL$ of sites and a configuration space $\Omega:=E^{\bL}$. 
In fact, we use the Dobrushin uniqueness condition in the two cases of $\bL$ which are $\bL\in\{-\N, \Z\}$. 

Consider a uniformly absolutely convergent (UAC) interaction $\Phi=\{\Phi_\Lambda(\cdot): \Lambda\Subset\bL\}$ on $\Omega$ and let $\gamma^{\Phi}$ be the corresponding Gibbsian specification. For any sites $i,j\in\bL$, define 
$$
C(\gamma^{\Phi})_{i,j}:=\sup_{\eta_{\bL\setminus\{j\}}=\overline{\eta}_{\bL\setminus\{j\}}}||\gamma_{\{i\}}^{\Phi}(\cdot|\eta)-\gamma^{\Phi}_{\{i\}}(\cdot|\overline{\eta})||_{\infty},
$$ 
where $||\cdot||_{\infty}$ is the supremum norm on $\mathcal{M}(\Omega)$ defined by $||\tau||_{\infty}:=\sup_{B\in\F}|\tau(B)|$  for any finite signed Borel measure $\tau$.
The infinite matrix  $C(\gamma^{\Phi}):=(C(\gamma^{\Phi})_{i,j})_{i,j\in\bL}$ is called the \emph{Dobrushin interdependence matrix}.
\medskip

\begin{defin}\cites{Georgii-book}
The specification $\gamma^{\Phi}$ satisfies the \textbf{\textit{Dobrushin uniqueness condition}} if
\begin{equation}\label{DUC for specification}
    c(\gamma^{\Phi}):=\sup_{i\in\bL}\sum_{j\in\bL}C(\gamma^{\Phi})_{i,j}<1. 
\end{equation}
\end{defin}

The Dobrushin uniqueness condition admits slightly stronger --- and easy to check--- forms.
One of them is a high-temperature condition, which 
is considered in \cite{EFMV2024}.
Another version is known as the strong magnetic field condition, which
we formulate here:
\begin{prop}\cite{Georgii-book}*{Example 8.13}\label{prop: DUC for strong magnetic fields} Let $\bL$ be any countable set and $E=\{\pm 1\}$.
Assume $\Phi$ be a UAC interaction on $\Omega=\{\pm 1\}^\bL$
such that $\Phi_{\{i\}}=-h \sigma_i$ for all $i\in\bL$ and for some $h\in \R$.
Suppose 
\begin{equation}\label{DUC for strong magnetic fields}
    |h|
    >
    \log \sup_{i\in\bL}
    \Big\{
\exp{\Big(
\frac{1}{2}\sum_{A\ni i, |A|\geq 2}\delta(\Phi_A)
\Big)}
\cdot
\sum_{A\ni i}(|A|-1)\delta(\Phi_A)
    \Big\},
\end{equation}
where  $\delta(f):=\sup\{|f(\xi)-f(\eta)|:\xi,\eta\in E^{\bL}\}$ is the variation of $f:E^{\bL}\to\R$. 
Then $\gamma^{\Phi}$ satisfies the Dobrushin uniqueness condition.

\end{prop}

The proof of Proposition \ref{DUC for strong magnetic fields} boils down to 
showing that for all $i,j$, $i\ne j$,
$$
C(\gamma^{\Phi})_{ij}
\leq
\exp{\Big(
-|h|+
\frac{1}{2} \sum_{A\ni i, |A|\geq 2}\delta(\Phi_A)
\Big)}
\cdot
\sum_{A\supset \{i,j\}} \delta(\Phi_A)
=:
\bar C(\Phi)_{ij},
$$
and hence
$c(\gamma^{\Phi})=\sup_i \sum_{j} C(\gamma^\Phi)_{ij}\le\sup_i \sum_{j} \bar C(\Phi)_{ij}=\bar c(\Phi)$.
Note that, here without loss of the generality, we set $\bar C(\Phi)_{ii}=0$ for all $i\in\bL$.
Notice that the non-negative matrix $\bar C(\Phi):=(\bar C(\Phi)_{i,j})_{i,j\in\bL}$  is symmetric.  

The condition (\ref{DUC for strong magnetic fields}) is stable under a \textit{perturbation} of the underlying model/interaction $\Phi$.
Indeed, let $\Psi=\{\Psi_V\}_{V\Subset\bL}$ be an interaction such that for $V\Subset\bL$, either $\Psi_V=\Phi_V$ or $\Psi_V=0$, then it is straightforward to check that the matrix $\bar C(\Psi)$ is dominated by $\bar C(\Phi)$, i.e., for all $i,j\in\bL$,
$\bar C(\Psi)_{ij}\leq \bar C(\Phi)_{ij}$, thus, in particular,
$\bar c(\Psi)\leq \bar c(\Phi)$. 
Hence $\Psi$ inherits the Dobrushin uniqueness condition from $\Phi$ as long as $\Phi$ satisfies (\ref{DUC for strong magnetic fields}).

The crucial property of the Dobrushin uniqueness condition is that it provides the uniqueness of the compatible probability measures with the specification $\gamma^\Phi$. 
In fact, we have the following theorem.
\begin{thm}\label{uniqueness under DU}\cite{Georgii-book}*{Chapter 8} If  $\gamma^{\Phi}$ satisfies the Dobrushin uniqueness condition (\ref{DUC for specification}), then
$\left|\G(\gamma^{\Phi})\right|= 1$.
\end{thm}

The validity of Dobrushin's criterion yields several important properties of the unique Gibbs state such as 
concentration inequalities and explicit bounds on the decay of correlations.
\\
The first property involves the coefficient $\bar{c} (\Phi)$, which governs the deviation behaviour of the unique Gibbs measure. 
Let
\begin{equation}
\delta_kF:=\sup\bigl\{F(\xi)-F(\eta):\xi_j=\eta_j,\; j\in\bL\setminus\{k\}\bigr\}
\end{equation}
denote the oscillation of a local function $F:\Omega\to\R$ at a site $k\in\bL$ and $\underline{\delta}(F)=(\delta_kF)_{k\in\bL}$ is the oscillation vector, where $\Omega=E^{\bL}$.

\begin{thm}\cites{K2003}\cite{V2018}\label{GCB and MCB}
    Suppose $\Phi$ is a UAC interaction satisfying (\ref{DUC for strong magnetic fields}) and let 
    $\mu_{\Phi}$ be its unique Gibbs measure. 
    Set 
    \begin{equation}
    D:=\frac{4}{(1-\bar{c}(\Phi))^2}\;.
    \end{equation}
   Then $\mu_{\Phi}$ satisfies a \textit{\textbf{Gaussian Concentration Bound}} with  the constant $D$, i.e., for any continuous function $F$ on $\Omega=E^{\bL}$, one has
    \begin{equation}\label{GCB ineq}
     \int_{\Omega}e^{F-\int_{\Omega}Fd\mu_\Phi}d\mu_\Phi\leq e^{D||\underline{\delta}(F)||^2_{2}}.
    \end{equation}
\end{thm}

The second important consequence of the Dobrushin's uniqueness condition is the following:
\begin{thm}\label{thm about comparing different specs}\cite{Georgii-book}*{Theorem 8.20}
Let $\gamma$ and $\tilde \gamma$ be two specifications on $\Omega=E^{\bL}$. 
Suppose $\gamma$ satisfies Dobrushin's condition.
For each $i\in\bL$, we let $b_i$ be a measurable function on $\Omega$ such that 
\begin{equation}
    ||\gamma_{\{i\}}(\cdot|\omega_{\{i\}^c})-\tilde\gamma_{\{i\}}(\cdot|\omega_{\{i\}^c})||_{\infty}\leq b_i(\omega)
\end{equation}
for all $\omega\in\Omega$.
If $\mu\in\G(\gamma)$ and $\tilde\mu\in\G(\tilde\gamma)$ then for all $f\in C(\Omega)$,
\begin{equation}
    |\mu(f)-\tilde\mu(f)|\leq \sum_{i,j\in\bL}\delta_i(f) D(\gamma)_{i,j}\tilde\mu(b_j).
\end{equation}
    
\end{thm}

We present the third property in the particular setup $\bL=\Z$, and it involves the $\Z\times\Z$ matrix
\begin{equation}
D(\gamma^{\Phi})=\sum_{n=0}^{\infty} C(\gamma^{\Phi})^n\;.
\end{equation}
 The sum of the  $\Z\times\Z$ matrices in the right-hand side converges due to the Dobrushin condition (\ref{DUC for specification}).

\begin{prop}\cites{Follmer1982, Georgii-book} Consider a UAC interaction $\Phi$ on $X=E^\Z$.
Assume the specification  $\gamma^{\Phi}$ satisfies the Dobrushin condition (\ref{DUC for specification}) and let $\mu$ be its unique Gibbs measure.
    Then, for all $f,g\in C(X)$ and $i\in\Z$, 
    \begin{equation}\label{Follmer's result for each term}
    \Big|\cov_{\mu}(f,g\circ S^i)\Big| 
    \leq\frac{1}{4}\sum_{k,j\in\Z} D(\gamma^{\Phi})_{jk}\cdot  \delta_kf\cdot\delta_{j-i}g.
    \end{equation}

\end{prop}
Suppose $\Phi$ satisfies the Dobrushin condition (\ref{DUC for strong magnetic fields}) and define the  non-negative symmetric $\Z\times \Z-$matrix by
    \begin{equation}\label{bar D for Phi}
    \bar D(\Phi):=\sum_{n\geq 0}\bar C(\Phi)^n\;.
    \end{equation}
The sum in the right-hand side of (\ref{bar D for Phi}) converges due to (\ref{DUC for strong magnetic fields}).
Furthermore, $\bar D(\Phi)$ is invertible and $\bar D(\Phi)^{-1}=I_{\Z\times \Z}-\bar C(\Phi)$, where $I_{\Z\times\Z}$ denotes the identity matrix (operator) on $l^2(\Z)$.

\noindent
In general, one can show that the off-diagonal elements $\bar D(\Phi)$ decay as they get far away from the diagonal, i.e., $\bar D(\Phi)_{ij}=o(1)$ as $|i-j|\to\infty$.
However, the decay rate largely depends on the interaction $\Phi$.
Now we discuss the decay rate in the case of Dyson interaction $\Phi$ defined by (\ref{eq: Dyson interaction}).
For $i\neq j$, one can readily check that
\begin{equation}
    \bar C(\Phi)_{ij}
    =
    \frac{2\beta\cdot \exp{(-|h|+2\beta \zeta (\alpha))}}{|i-j|^\alpha}=O(|i-j|^{-\alpha}).
\end{equation}
Thus by Jaffard's theorem (\cite{Jaffard1990}*{Proposition 3} and \cite{GK2010}*{Theorem 1.1}), one has the same asymptotics for the off-diagonal elements of the inverse matrix $(I_{\Z\times\Z}-\bar C(\Phi))^{-1}=\bar D(\Phi)$,
i.e., $\bar D(\Phi)_{ij}=O(|i-j|^{-\alpha})$.
Thus there exists a constant $C_J>0$ such that for all $i,j\in\Z$,
    \begin{equation}\label{eq: result from Jaffard}
        \bar D(\Phi)_{i,j}\leq \frac{C_J}{(1+|i-j|)^\alpha}.
    \end{equation}

\section{Construction of intermediate interactions}\label{Section: Construction of intermediate interactions}
In this section, we recall the construction of the intermediate interaction from \cite{EFMV2024}.
For the Dyson interaction $\Phi$ given by (\ref{eq: Dyson interaction}) and each $k\in\Z_+$, we construct intermediate interactions ${\Psi}^{(k)}$ in the following way.
We will represent $\Z=-\mathbb N\cup \Z_+$.
Consider the countable collection of finite subset of $\Z$ with  endpoints in $-\mathbb N$ and $\Z_+$:
$$
\mathcal A=\{\Lambda\Subset \Z:\, \min(\Lambda)<0,\ \max \Lambda\ge 0\}.
$$
Index elements of $\mathcal A$ in an order $\mathcal A=\{\Lambda_1,\Lambda_2,\ldots\}$ in a way that for every $N\in\N$, there exists $k_N\in\N$ such that
\begin{equation}
\sum_{i=1}^{k_N}\Phi_{\Lambda_i}
=
\sum_{\substack{
\min V<0\leq \max V\\
V\subset [-N,N]}} \Phi_V.       
\end{equation}
Then define
$$
{\Psi}^{(k)}_\Lambda = \begin{cases}  \Phi_{\Lambda},&\quad \Lambda\notin\{\Lambda_i: i\geq k+1\},\\
0, &\quad \Lambda\in\{\Lambda_i: i\geq k+1\},
\end{cases}\;
$$
In other words, we first remove all $\Phi_\Lambda$'s with $\Lambda\in\mathcal A$ from $\Phi$, and then add them one by one. 
Clearly, all the constructed interactions are UAC. 

\begin{remark}\label{convergence of intermediate interactions and specs }
    \begin{itemize}
        \item[1)] Every $\Psi^{(k)}$ is a local (finite) perturbation of $\Psi^{(0)}$, and $\Psi^{(k)}$ tend to $\Phi$
        as $k\to\infty$, in the sense that 
        $\Psi^{(k)}_{\Lambda}\rightrightarrows\Phi_{\Lambda}$ 
        for all $\Lambda\Subset\Z$.
        \item[2)] For any finite volume $V$, one can readily check that 
        $$||H^{\Psi^{(k)}}_V-H^{\Phi}_V||_{\infty}\leq \sum\limits_{\substack{\Lambda_j\cap V\neq \emptyset \\j\geq k}}||\Phi_{\Lambda}||_{\infty}\xrightarrow[k\to\infty]{}0. 
        $$
        \item[3)] For specifications, it can also be concluded that $\gamma^{\Psi^{(k)}}$ converges to $\gamma^{\Phi}$ as $k\to\infty$, i.e., for all $B\in\F$ and $V\Subset\Z$,
        $$
\gamma^{\Psi^{(k)}}_V(B|\omega)\xrightarrow[k\to\infty] 
    {}\gamma^{\Phi}_V(B|\omega)\; \text{ 
        uniformly on the b.c. } 
        \omega\in\Omega.
        $$
        \item[4)] In addition, if $\nu^{(k)}$ is a Gibbs measure for $\Psi^{(k)}$, then by Theorem 4.17 in \cite{Georgii-book}, any weak$^*$-limit point of the sequence $\{\nu^{(k)}\}_{k\geq 0}$ becomes a Gibbs measure for the potential $\Phi$.

    \end{itemize}
\end{remark}

Another important observation is the following: since we have constructed $\Psi^{(0)}$ from $\Phi$ by removing all the interactions between the left $-\mathbb N$ and the right $\Z_+$ half-lines, the corresponding specification $\gamma^{\Psi^{(0)}}$ becomes product type \cite{Georgii-book}*{Example 7.18}. 
More precisely, $\gamma^{\Psi^{(0)}}=\gamma^{\Phi^-}\times\gamma^{\Phi^+}$, where $\Phi^-$ and $\Phi^+$ are the restrictions of $\Phi$ to the half-lines $-\N$ and $\Z_+$, respectively. Thus we have the following for the extreme Gibbs measures \cite{Georgii-book}*{Example 7.18}:

\begin{equation}\label{extreme Gibbs measures for Psi^0}
    \text{ex}\;\G(\gamma^{\Psi^{(0)}})
    =\{\nu^-\times\nu^+: \nu^-\in \text{ex}\;\G(\gamma^{\Phi^-}),\;\; \nu^+\in \text{ex}\;\G(\gamma^{\Phi^+}) \}
\end{equation}

For $k\geq 0$, let $\nu^{(k)}\in\G(\Psi^{(k)})$. 
For any $k\geq 1$, we also consider the following function:
\begin{equation}\label{densities with respect to the product of one-sided Gibbs measures}
    f^{(k)}=\frac{e^{-\sum_{i=1}^k\Phi_{\Lambda_i}}}{\int_{\Omega}e^{-\sum_{i=1}^k\Phi_{\Lambda_i}}d\nu^{(0)}}.
\end{equation}

Now we state two important lemmas which will be used in the proof of Theorem \ref{main result}.

\begin{lem}\label{lemma: equiv of intermediate Gibbs meas}\cite{EFMV2024}*{Theorem 6.5}
If the interaction $\Phi$ satisfies the Dobrushin uniqueness condition (\ref{DUC for strong magnetic fields}), so do all interactions $\Psi^{(k)}$. 
Furthermore, $\nu^{(k)}$ and $\nu^{(0)}$ are equivalent with
    \begin{equation}\label{densities with respect to the product of half-line Gibbs measures}
   \frac{d\nu^{(k)}}{d\nu^{(0)}}=\frac{e^{-\sum_{i=1}^k\Phi_{\Lambda_i}}}{\int_{X}e^{-\sum_{i=1}^k\Phi_{\Lambda_i}}d\nu^{(0)}} 
   \quad,\quad   \frac{d\nu^{(0)}}{d\nu^{(k)}}=\frac{e^{\sum_{i=1}^k\Phi_{\Lambda_i}}}{\int_{X}e^{\sum_{i=1}^k\Phi_{\Lambda_i}}d\nu^{(k)}}\;.
\end{equation}
\end{lem}

Lemma \ref{lemma: equiv of intermediate Gibbs meas} further implies that for any $k\geq 1$,
\begin{equation}\label{densities-the recursion}
    \frac{d\nu^{(k)}}{d\nu^{(k-1)}} = \frac {e^{-\Phi_{\Lambda_k}}}
 { \int e^{-\Phi_{\Lambda_{k}}} d\nu^{(k-1)}
 }.
\end{equation}

\begin{lem}\label{lemma: UI and Dunford-Pettis}\cite{EFMV2024}*{Theorem D}
    Assume that $\Phi$ satisfies the Dobrushin uniqueness condition (\ref{DUC for strong magnetic fields}). 
    Suppose the family $\{f^{(k)}\}_{k\in\N}$ is uniformly integrable in $L^1(\nu^{(0)})$. 
    Then the weak$^*$ limit point of the sequence $\{\nu^{(k)}\}$ is a Gibbs measure for $\Phi$ and absolutely continuous with respect to $\nu^{(0)}$.
\end{lem}

\noindent
Proofs of both Lemma \ref{lemma: equiv of intermediate Gibbs meas} and Lemma \ref{lemma: UI and Dunford-Pettis} can be found in \cite{EFMV2024}.

We end this section by noting that if $h\geq 0$, then the Griffiths–Kelly–Sherman (GKS) inequality applies to the Gibbs measures $\nu^{(k)}$ and $\mu$ \cites{Griffiths1966, KS1967, Griffiths1968, Ginibre1970}.
Namely, for every $\tau\in\{\mu, \nu^{(k)}, \; k\geq 0\}$ and all  
$A, B\Subset \Z$, one has
\begin{equation}\label{eq: GKS-1}
\int_{\Omega}\sigma_Ad\tau\geq 0
\end{equation}
and 
\begin{equation}\label{eq: GKS-2}
    \int_{\Omega}\sigma_A\sigma_Bd\tau\geq \int_{\Omega}\sigma_Ad\tau\cdot\int_{\Omega}\sigma_Bd\tau,
\end{equation}
where $\sigma_A:=\prod_{i\in A}\sigma_i$ and for $i\in \Z$, $\omega\in X\mapsto \sigma_i(\omega)=\omega_i$.
The GKS inequalities also imply the following inequality for $\tau\in\{\mu, \nu^{(k)}, \; k\geq 0\}$: for every $A, B\Subset\Z$ and $t\geq 0$,
\begin{equation}\label{eq: Griffith important cor}
    \int_{\Omega} \sigma_B e^{t\sigma_A}d\tau\geq \int_{\Omega}\sigma_Bd\tau \cdot\int_{\Omega}e^{t\sigma_A}d\tau.
\end{equation}

\section{Proof of Theorem \ref{main result}}\label{Section: Proof of main result}
\textbf{Part (i):}
Firstly, note that the uniqueness of Gibbs measures for $\Phi$ and $\Psi^{(0)}$ easily follows from Theorem \ref{uniqueness under DU}.
Thus by (\ref{extreme Gibbs measures for Psi^0}), the restriction $\Phi^+$ of $\Phi$ to the half-line $\Z_+$ also has a unique Gibbs measure $\nu^+$.
Furthermore, Theorem 4.8 in \cite{CLS2019} implies that the transfer operator $\La_\phi$ has a unique eigenprobability $\nu$, and $\nu=\nu^+$ (see also Subsection 2.2 in \cite{EFMV2024}).

\bigskip

\textbf{Part (ii):} 
We start the proof of the second part by noting that the positive and negative fields are related by the global-spin flip transformation $\s: X_+\to X_+$ defined by $\s(x)_n=-x_n$ for all $n\in\Z_+$. 
In fact, for the Dyson potential $\phi$ given by (\ref{eq: Dyson potential on half-line}), one has that 
\begin{equation}
    \phi\circ\s=-h\sigma_0 + \beta\sum_{n=1}^\infty\frac{\sigma_0\sigma_n}{n^\alpha}.
\end{equation}
For any continuous potential $\psi\in C(X_+)$ one can easily check that
\begin{equation}\label{identity ito transfer op for spin-flip transf}
    \big(
    \La_{\psi\circ\s} \;\mathfrak{f}
    \big)
    \circ \s
    =
    \La_\psi(\mathfrak{f}\circ\s),\;\; \mathfrak{f}\in \R^{X_+}.
\end{equation}
Thus, the topological pressures of $\psi$ and $\psi\circ\s$ coincide, i.e., $p(\psi)=p(\psi\circ \s)$ and $\nu\in\M_{1}(X_+)$ is an eigenprobability of $\La_\psi^*$ if and only if $\nu\circ \s$ is an eigenprobability of $\La_{\psi\circ \s}^*$.
Furthermore, $\mathtt{h}$ is an eigenfunction of $\La_\psi$ corresponding to $\lambda=e^{p(\psi)}$ if and only if $\h\circ \s$ is an eigenfunction of $\La_{\psi\circ \s}$ corresponding to $\lambda$.

By the above argument, we may assume without loss of generality that the external field $h$ is positive for the remainder of the proof.

We below show that the unique Gibbs measure $\mu$ for $\Phi$ is absolutely continuous with respect to $\nu^{(0)}$.
Then this yields that the restriction $\mu_+$ of $\mu$ to $\Z_+$ is absolutely continuous with respect to $\nu$ and the Radon-Nikdoym density $f_+:=\frac{d\mu_+}{d\nu}$ is given by
\begin{equation}
    f_+(\omega_{0}^\infty)=\int_{X_-} f(\xi_{-\infty}^{-1}\omega_{0}^\infty) \nu_-(d\xi),
\end{equation}
where $f=\frac{d\mu}{d\nu^{(0)}}$.
Then Proposition 3.1 in \cite{EFMV2024} also yields that the Radon-Nikodym density $\frac{d\mu_+}{d\nu}$ is an eigenfunction of the transfer operator $\La_\phi$ corresponding to $\lambda=e^{P(\phi)}$.

Note that for the Dyson interaction $\Phi$, one has, for all $i\in\N,j\in\Z_+$ and $k\in\Z$, 
\begin{equation*}
    \delta_k(\Phi_{\{-i,j\}})=\begin{cases} 
      0, & k\notin \{-i,j\}; \\
      \frac{2\beta}{(i+j)^{\alpha}}, & k\in\{-i,j\}.
     \end{cases}
    \end{equation*}
Thus for all $i\in\N,\; j\in\Z_+$, $||\underline\delta(\Phi_{\{-i,j\}})||^2_2=\frac{8\beta^2}{(i+j)^{2\alpha}}$. 
Therefore, $\sum_{\substack{i\in\N,\\j\in\Z_+}}||\underline\delta(\Phi_{\{-i,j\}})||^2_2
=
\sum_{\substack{i\in\N,\\j\in\Z_+}}\frac{8\beta^2}{(i+j)^{2\alpha}}<\infty$.

Consider any $\Lambda_k\in\A$ and $n\in\N$. It follows from (\ref{densities-the recursion}) that
$$
\int_{\Omega} \Phi_{\Lambda_k}d\nu^{(n)}=
\frac{-J_{\Lambda_k}\int_{\Omega}\sigma_{\Lambda_k}e^{J_{\Lambda_n}\sigma_{\Lambda_n}}d\nu^{(n-1)}}{\int_{\Omega}e^{J_{\Lambda_n}\sigma_{\Lambda_n}}d\nu^{(n-1)}},
$$
here for $\Lambda\Subset\Z_+$, $J_\Lambda:=h$ if $\Lambda=\{i\}$ and $J_\Lambda:=\frac{\beta}{|i-j|^\alpha}$ if $\Lambda=\{i,j\}$.
Then by (\ref{eq: Griffith important cor}),
\begin{equation}\label{monotonicity of the means nu tilde}
\int_{\Omega} \Phi_{\Lambda_k}d\nu^{(n)}\leq 
\frac{\int_{\Omega}\Phi_{\Lambda_k}d\nu^{(n-1)}\cdot \int_{\Omega}e^{J_{\Lambda_n}\sigma_{\Lambda_n}}d\nu^{(n-1)}}{\int_{\Omega}e^{J_{\Lambda_n}\sigma_{\Lambda_n}}d\nu^{(n-1)}}=
\int_{\Omega}\Phi_{\Lambda_k}d\nu^{(n-1)}.
\end{equation}
Note that $\nu^{(n)}$ converges to $\mu$ as $n\to\infty$ in the weak star topology. Thus one can obtain from (\ref{monotonicity of the means nu tilde}) that for all $\Lambda_k\in\A$, $n\in\N$,
\begin{equation}\label{ub lb bounds for nu tilde mean}
    \int_{\Omega} \Phi_{\Lambda_k}d\nu^{(0)}\geq \int_{\Omega} \Phi_{\Lambda_k}d\nu^{(n)}
    \geq\lim_{n\to\infty}\int_{\Omega} \Phi_{\Lambda_k}d\nu^{(n)}
    = \int_{\Omega} \Phi_{\Lambda_k}d\mu.
\end{equation}

For all $k\in\N$, denote $W_k:=\sum_{i=1}^k\Phi_{\Lambda_i}$. 
Then we aim to prove 
\begin{equation}\label{place to apply Valle-Poissen}
    \sup_{k\geq 0} \int_{X}f^{(k)}\log f^{(k)}d\nu^{(0)}<\infty.\;
\end{equation}
Hence by applying de la Vall\'ee Poussin's theorem to the family $\{f^{(k)}:k\in\N\}$ and to
the function $t\in(0,+\infty)\mapsto t\log t$, one concludes that the family  $\{f^{(k)}:k\in\N\}$ is uniformly integrable in $L^1(\nu^{(0)})$. 
Then Theorem \ref{lemma: UI and Dunford-Pettis} imply that the unique limit point $\mu$ is absolutely continuous with respect to $\nu^{(0)}$.

Note that  for every $k\in\N$,
\begin{equation}\label{the sum}
    \int_{X}f^{(k)}\log f^{(k)}d\nu^{(0)}
    =
    \int_{X}-W_kd\nu^{(k)}
    -
    \log \int_{X}e^{-W_k}d\nu^{(0)}.
\end{equation} 
Here one can easily check that 
$$
\sup_{k\in\N}\int_{X}e^{-W_k}d\nu^{(0)}
=
+\infty.
$$
Under conditions of our main theorem, we argue that 
$$
\sup_{k\in\N}\int_{X}-W_kd\nu^{(k)}
=
+\infty.
$$
Therefore, the method of the proof of Theorem E in \cite{EFMV2024} does not work here.
We instead prove the following two claims:
\medskip

\noindent
{\bf Claim 1:}
\[
\sup_{k\in\N}\Big|\int_{X}W_kd\nu^{(k)}-\int_X W_k d\mu\Big|<\infty;
\]
\noindent
{\bf Claim 2:}
\begin{equation*}\label{UI of numerator}
    \sup_{k\in\N}\Big| 
    -\log \int_{X}e^{-W_k}d\nu^{(0)} 
    -
    \int_X W_kd\mu
    \Big|
    <\infty.
\end{equation*}
One can easily see from (\ref{the sum}) that Claim 1 and Claim 2 indeed imply (\ref{place to apply Valle-Poissen}).

Note that 
 for all $k\in\N$,  $\bar c (\Psi^{(0)})\leq \bar c(\Psi^{(k)})\leq \bar c(\Phi)$.
Therefore, the Dobrushin uniqueness condition $\bar c(\Phi)<1$ is inherited by all the intermediate interactions.  Applying the first part of Theorem \ref{GCB and MCB} we see that the (only) measure 
$ \mu\in\G(\Phi)$ and all the intermediate measures $\nu^{(k)}$, $k\geq 0$, satisfy the Gaussian Concentration Bound with the same constant 
$D:=\frac{4}{(1-\bar c(\Phi))^2}$. 
This implies that, for all $k\in\N$, 
\begin{equation}
    \int_{X}e^{-\Phi_{\Lambda_k}}d\nu^{(k-1)}\leq e^{ D||\underline\delta(\Phi_{\Lambda_k})||^2_2} e^{-\int_{X}\Phi_{\Lambda_k}d\nu^{(k-1)}}.
\end{equation}
We combine this inequality with (\ref{densities-the recursion}) to iterate
\begin{eqnarray}
\lefteqn{ \int_{X}e^{-(\Phi_{\Lambda_k}+\Phi_{\Lambda_{k-1}})}d\nu^{(k-2))}}\nonumber\\
&=&
  \int_{X}e^{-\Phi_{\Lambda_{k}}}d\nu^{(k-1))}\, \int_{X}e^{-\Phi_{\Lambda_{k-1}}}d\nu^{(k-2))} \nonumber\\[8pt]
 &\le& 
 e^{ D(||\underline\delta(\Phi_{\Lambda_k})||^2_2+||\underline\delta(\Phi_{\Lambda_{k-1}})||^2_2)} \cdot e^{-(\int_{X}\Phi_{\Lambda_k}d\nu^{(k-1)}+\int_{X}\Phi_{\Lambda_{k-1}}d\nu^{(k-2)})}.
\end{eqnarray}
By induction this yields 
\begin{equation}
    \int_{X}e^{-\sum_{i=1}^k\Phi_{\Lambda_i}}d\nu^{(0)}\leq 
    e^{ D\sum_{i=1}^k||\underline\delta(\Phi_{\Lambda_i})||^2_2}\cdot e^{-\sum_{i=1}^k\int_{X}\Phi_{\Lambda_i}d\nu^{(i-1)}}.
\end{equation}
Similarly, one can also obtain the lower bound:
\begin{equation}
    \int_{X}e^{-\sum_{i=1}^k\Phi_{\Lambda_i}}d\nu^{(0)}
    \geq
    e^{- D\sum_{i=1}^k||\underline\delta(\Phi_{\Lambda_i})||^2_2}\cdot e^{-\sum_{i=1}^k\int_{X}\Phi_{\Lambda_i}d\nu^{(i-1)}}.
\end{equation}
Hence for all $k\in \N$,
\begin{equation}
    -C_1- \sum_{i=1}^k\int_X \Phi_{\Lambda_i}d\nu^{(i-1)}
    \leq
    \log\int_{X} e^{-\sum_{i=1}^k\Phi_{\Lambda_i}}d\nu^{(0)}
    \leq
    C_1- \sum_{i=1}^k\int_X \Phi_{\Lambda_i}d\nu^{(i-1)}, 
\end{equation}
where
$
C_1
:= D\cdot\sum_{i=1}^\infty||\underline\delta(\Phi_{\Lambda_i})||^2_2
=
\sum_{\substack{i\in\N,\\j\in\Z_+}}\frac{8\beta^2}{(i+j)^{2\alpha}}<\infty.
$
Thus instead of Claim 2, it suffices to show the following:
\begin{equation}\label{eq: alternative for the denominator}
    \sup_{k\in\N}
    \Big| 
    \sum_{i=1}^k\int_X \Phi_{\Lambda_i}d\nu^{(i-1)} 
    -
    \int_X W_kd\mu
    \Big|
    <\infty.
\end{equation}
Thus, in the light of (\ref{ub lb bounds for nu tilde mean}), Claim 1 and Claim 2 are implied by the following:

\noindent
{\bf Claim 3:}
\begin{equation}\label{eq: nu^0 and mu left-right sums comparison}
    \sup_{k\in\N}
    \Big| 
    \int_X W_kd\nu^{(0)} 
    -
    \int_X W_kd\mu
    \Big|
    <\infty.
\end{equation}
By the GKS inequality (c.f., \ref{ub lb bounds for nu tilde mean}), we, in fact, have that 
$
\int_X \Phi_{\Lambda_i}d\nu^{(0)} 
    -
\int_X \Phi_{\Lambda_i}d\mu\geq 0,
$
hence the sequence 
$\{\int_X W_kd\nu^{(0)} 
    -
\int_X W_kd\mu\}_k$ increases with $k$.
Thus instead of considering all $W_k$'s, it is enough to consider the subsequence $\{W_{[N]}\}_{N\in\N}$s, where 
$$
W_{[N]}(\omega)
:=
\sum_{i=1}^{N}\sum_{j=0}^{N} \Phi_{\{-i,j\}}(\omega)
=
\sum_{i=1}^{N}\sum_{j=0}^{N}-\frac{\beta \omega_{-i}\omega_{j}}{(i+j)^{\alpha}}.
$$
Then
\begin{eqnarray}
    \int_X W_{[N]}d\nu^{(0)} 
    -
    \int_X W_{[N]}d\mu
    \notag&=&
    \sum_{i=1}^{N}\sum_{j=0}^{N}\int_X -\frac{\beta \omega_{-i}\omega_{j}}{(i+j)^{\alpha}}\nu^{(0)}(d\omega)
    -
    \sum_{i=1}^{N}\sum_{j=0}^{N}\int_X-\frac{\beta \omega_{-i}\omega_{j}}{(i+j)^{\alpha}}\mu(d\omega)\\
    \notag &=&
    \sum_{i=1}^{N}\sum_{j=0}^{N} \frac{-\beta\nu_-(\sigma_{-i})\cdot\nu_+(\sigma_j)}{(i+j)^\alpha}\\
    \notag&-&
    \sum_{i=1}^{N}\sum_{j=0}^{N}\frac{-\beta\mathtt{cov}_\mu(\sigma_{-i},\sigma_j)-\beta\mu(\sigma_0)^2}{(i+j)^\alpha}\\
    &=&\label{sum of correlations}
    \sum_{i=1}^{N}\sum_{j=0}^{N}\frac{\beta\mathtt{cov}_\mu(\sigma_{-i},\sigma_j)}{(i+j)^\alpha}\\
    &+&\label{problematic double sum}
    \beta\sum_{i=1}^{N}\sum_{j=0}^{N}\frac{\mu(\sigma_0)^2-\nu_-(\sigma_{-i})\cdot\nu_+(\sigma_j)}{(i+j)^\alpha},
\end{eqnarray}
here we used the fact that $\nu^{0}=\nu_-\times\nu_+$.
By (\ref{Follmer's result for each term}) and (\ref{eq: result from Jaffard}), we know that in the Dobrushin uniqueness region (\ref{sum of correlations}) is bounded, i.e., 
$$
\sum_{i=1}^{\infty}\sum_{j=0}^{\infty}\frac{\beta|\mathtt{cov}_\mu(\sigma_{-i},\sigma_j)|}{(i+j)^\alpha}<\infty.
$$
Therefore, to show 
\begin{equation}
    \sup_{N\in\N}
    \Big(
    \int_X W_{[N]}d\nu^{(0)} 
    -
    \int_X W_{[N]}d\mu
    \Big)
    <
    \infty
\end{equation}
it is enough to show that
\begin{equation}
    \beta\sum_{i=1}^{\infty}\sum_{j=0}^{\infty}\frac{|\mu(\sigma_0)^2-\nu_-(\sigma_{-i})\cdot\nu_+(\sigma_j)|}{(i+j)^\alpha}
    <
    \infty.
\end{equation}

By Theorem A in \cite{EFMV2024}, we have that $\nu_+(\sigma_j)\xrightarrow[j\to\infty]{}\mu(\sigma_0)$, similarly, one can also argue that 
$\nu_-(\sigma_{-i})\xrightarrow[i\to\infty]{}\mu(\sigma_0)$.
For $n\in\Z$, denote 
\begin{equation}
    t_n:=|\nu^{(0)}(\sigma_n)-\mu(\sigma_0|.
\end{equation}
Then (\ref{problematic double sum}) can be rewritten in the following form
\begin{equation}\label{imp sum to be estimated}
\beta\sum_{i=1}^{N}\sum_{j=0}^{N}\frac{|\mu(\sigma_0)^2-\nu_-(\sigma_{-i})\cdot\nu_+(\sigma_j)|}{(i+j)^\alpha}
=
\beta\sum_{i=1}^{N}\sum_{j=0}^{N}\frac{|\mu(\sigma_0) t_j+\mu(\sigma_0)t_{-i}-t_{-i}t_j|}{(i+j)^\alpha}.
\end{equation}

Now we show that there exists a constant $C_2=C_2(\alpha,\beta, h)>0$ such that for every $n\in\Z$,
\begin{equation}\label{eq: speed of convergence of t_n}
    t_n\leq \frac{C_2}{(|n|+1)^{\alpha-1}}.
\end{equation}
In fact, we prove the above inequality for positive $n$'s and the case of negative $n$'s can be treated similarly.
Fix $j\in\N$.
Then by Theorem \ref{thm about comparing different specs}, one has that 
\begin{equation}\label{eq: est for t_j in terms of D}
    t_j
    =
    \Big|\int_X \sigma_jd\mu-\int_X \sigma_j d\nu^{(0)} \Big|
    \leq
    2(\bar D b)_j,
\end{equation}
where $\bar D=\bar D(\Phi)$ which is given by (\ref{bar D for Phi}) and $b=(b_s)_{s\in\Z}$ is given by the following:
\begin{equation}
    b_s:=\sup_{\omega\in X} 
    \lVert 
    \gamma_{\{s\}}^\Phi(\cdot|\omega_{\{s\}^c})-\gamma^{\Psi^{(0)}}_{\{s\}}(\cdot| \omega_{\{s\}^c}) 
    \rVert_{\infty}.
\end{equation}
For any $\xi, \omega\in X$, one has that 
\begin{multline*}
    \gamma_{\{s\}}^\Phi(\xi_s|\omega_{\{s\}^c})-\gamma^{\Psi^{(0)}}_{\{s\}}(\xi_s| \omega_{\{s\}^c})\\
    =
    \frac{e^{-H_s^\Phi(\xi_s\omega_{\{s\}^c})-H_s^{\Psi^{(0)}}(\xi_s\omega_{\{s\}^c})}
    \sum\limits_{\eta_s}
    \Big[
    e^{H_s^{\Psi^{(0)}}(\xi_s\omega_{\{s\}^c})-H_s^{\Psi^{(0)}}(\eta_s\omega_{\{s\}^c})}
    -
    e^{H_s^{\Phi}(\xi_s\omega_{\{s\}^c})-H_s^{\Phi}(\eta_s\omega_{\{s\}^c}) }
    \Big] }
    {\Big(
    \sum\limits_{\eta_s}e^{-H_s^{\Phi}(\eta_s\omega_{\{s\}^c})}
    \Big)
    \cdot
    \Big(
    \sum\limits_{\eta_s}e^{-H_s^{\Phi^{(0)}}(\eta_s\omega_{\{s\}^c})}
    \Big)}.
\end{multline*}
Thus since 
\begin{equation*}
   \sup_{\xi,\eta,\omega}\Big|
H_s^{\Phi}(\xi_s\omega_{\{s\}^c})-H_s^{\Phi}(\eta_s\omega_{\{s\}^c}) 
-
\Big(
H_s^{\Psi^{(0)}}(\xi_s\omega_{\{s\}^c})-H_s^{\Psi^{(0)}}(\eta_s\omega_{\{s\}^c})
\Big)
\Big|
=
O\Big(\frac{1}{(|s|+1)^{\alpha-1}}\Big) 
\end{equation*}
\begin{equation}
    \sup_{\omega\in X}||\gamma_{\{s\}}^\Phi(\cdot|\omega_{\{s\}^c})-\gamma^{\Psi^{(0)}}_{\{s\}}(\cdot| \omega_{\{s\}^c})||_\infty
    =
    O\Big(\frac{1}{(|s|+1)^{\alpha-1}}\Big).
\end{equation}
Hence there exists $C_3>0$ such that for all $s\in\Z$,
\begin{equation}
    b_s
    \leq
    \frac{C_3}{(|s|+1)^{\alpha-1}}.
\end{equation}
Then (\ref{eq: est for t_j in terms of D}) yields that 
\begin{equation}
    t_j
    \leq
    2C_3
    \sum_{s\in\Z} \frac{\bar D_{j,s}}{(|s|+1)^{\alpha-1}}.
\end{equation}
Hence by using (\ref{eq: result from Jaffard}), we obtain that 
\begin{equation}
    t_j
    \leq
    2C_3 C_J
    \sum_{s\in\Z} \frac{1}{(1+|s-j|)^\alpha(|s|+1)^{\alpha-1}}.
\end{equation}
One can see that 
\begin{equation}
    \sum_{s\in\Z} \frac{1}{(1+|s-j|)^\alpha(|s|+1)^{\alpha-1}}
    =
    O(j^{1-\alpha}). 
\end{equation}
Hence we conclude (\ref{eq: speed of convergence of t_n}) in the case of $n>0$.

Using (\ref{eq: speed of convergence of t_n}), we can now estimate (\ref{imp sum to be estimated}) from the above.
In fact, one has
\begin{multline*}
\sum_{i=1}^{N}\sum_{j=0}^{N}\frac{|\mu(\sigma_0) t_j+\mu(\sigma_0)t_{-i}-t_{-i}t_j|}{(i+j)^\alpha}\\
\leq 
\sum_{i=1}^{N}\sum_{j=0}^{N}
\frac
{|\mu(\sigma_0)| t_j}
{(i+j)^\alpha}
+
\sum_{i=1}^{N}\sum_{j=0}^{N}\frac{|\mu(\sigma_0)|t_{-i}}{(i+j)^\alpha}
+
\sum_{i=1}^{N}\sum_{j=0}^{N}\frac{t_{-i}t_j}{(i+j)^\alpha}.    
\end{multline*}
Now we only show that the first term remains bounded as $N\to\infty$, and the boundedness of the other terms can be shown similarly.
In fact, there exists $C_4>0$ such that for all $N$,
\begin{equation*}
\sum_{j=0}^{N}\sum_{i=1}^{N}
\frac
{t_j}
{(i+j)^\alpha} 
\leq
\sum_{j=0}^{N}\sum_{i=1}^{\infty}
\frac
{t_j}
{(i+j)^\alpha} 
\leq
C_4\sum_{j=1}^{N}
\frac{t_j}{j^{\alpha-1}}.
\end{equation*}
Thus by (\ref{eq: speed of convergence of t_n}) and by taking into account the fact that $\alpha>\frac{3}{2}$, one immediately concludes that for all $N\in\N$,
\begin{equation}
\sum_{j=0}^{N}\sum_{i=1}^{N}
\frac
{t_j}
{(i+j)^\alpha} 
\leq
C_2C_4\sum_{j=1}^N\frac{1}{j^{2\alpha-2}}
\leq
C_2C_4\sum_{j=1}^\infty\frac{1}{j^{2\alpha-2}}
<
\infty.
\end{equation}

\bigskip

\textbf{Part (iii):} By the first part, we can also conclude that the entire sequence $\{f^{(k)}\}_{k\in\N}$ converges to $f=\frac{d\mu}{d\nu^{(0)}}$ as $k\to\infty$ in the weak topology in $L^1(\nu^{(0)})$.
This follows from the fact that a limit point of the sequence $\{f^{(k)}\}_{k\in\N}$ should be a Radon-Nikodym density of some $\mu^*\in\G(\Phi)$, but $\G(\Phi)=\{\mu\}$, therefore, $\frac{d\mu}{d\nu^{(0)}}$ is the only limit point.
Thus 
\begin{equation}\label{eq: limit for two-sided RN-density}
    f=\lim_{k\to\infty}\frac{e^{-W_k}}{\int_{X}e^{-W_k}d\nu^{(0)}},
\end{equation}
here the limit is understood in the weak topology in $L^1(\nu^{(0)})$.

Note that the weak convergence of $\{f^{(k)}\}$ to $f$ in the weak topology in $L^1(\nu^{(0)})$ also implies the weak convergence of $f^{(k)}_+:=\int_{X_-}f^{(k)}d\nu_-$ to $f_+$ in the weak topology in $L^1(\nu)$.
In fact, for any bounded $g\in \F_{\Z_+}$, by taking into account the fact that $\nu^{(0)}=\nu_-\times\nu$, one has 
\begin{eqnarray*}
    \int_{X_+} (f^{(k)}_+-f_+)g d\nu
    &=&
    \int_{X_+}\Big[\int_{X_-}(f^{(k)}-f)d\nu_-\Big]g d\nu\\
    &=&
    \int_{X_+}\int_{X_-}g(f^{(k)}-f)d\nu_- d\nu\\
    &=&
    \int_{X}g(f^{(k)}-f)d\nu^{(0)}\xrightarrow[k\to\infty]{}0.
\end{eqnarray*}
Hence we obtain an analog of (\ref{eq: limit for two-sided RN-density}) for $f_+$, i.e.,

\begin{equation}\label{eq: limit for one-sided RN-density}
    f_+=\lim_{k\to\infty}\frac{\int_{X_-}e^{-W_k}d\nu_-}
    {\int_{X}e^{-W_k}d\nu^{(0)}}.
\end{equation}
Here the above limit should be again understood in an appropriate topology.

Now assume the transfer operator $\La_\phi$ has a continuous eigenfunction $\tilde f$, or in other words, the Radon-Nikodym derivative $\frac{d\mu_+}{d\nu}$ has a continuous version $\tilde f$.
Then for $\nu-$a.e. $x\in X_+$, $\tilde f(x)=f_+(x)$, and $\tilde f$, as a function on compact space $X_+$, is bounded.
Hence $f_+\in L^\infty(\nu)$,
therefore, 
\begin{equation}\label{imp sup to prove discontinuity}
    \sup\Big\{\frac{1}{\nu([y_\Lambda])}\int_{[y_\Lambda]} f_+\; d\nu: y\in X_+,\; \Lambda\Subset\Z_+ \Big\}
    \leq
    \lVert  f_+\rVert_{L^\infty(\nu)}
    <\infty.
\end{equation}
Hence to prove the second part of Theorem \ref{main result}, it is enough to show that the supremum in (\ref{imp sup to prove discontinuity}) is infinite.
In fact, below we show that 
\begin{equation}
    \sup_{n\in\N}\frac{1}{\nu([1_0^n])}\int_{[1_0^n]} f_+ d\nu
    =
    +\infty.
\end{equation}
In other words, we show that $1_{\Z_+}=(+1_i)_{i\in\Z_+}\in X_+$ is an \textit{essential discontinuity point} of the Radon-Nikdoym density $f_+=\frac{d\mu_+}{dd\nu}$.
To do so, fix $n\in \N$ and consider $g=\mathds{1}_{[1_0^n]}$.
By (\ref{eq: limit for one-sided RN-density}), one has that 
\begin{equation}\label{eq with all k f_+ 1_{[1_0^n]}}
    \int_{X_+} f_+ \mathds{1}_{[1_0^n]} d\nu =\lim_{k\to\infty}
    \int_{[1_0^n]}
    \frac
    {\int_{X_-}e^{-W_k(\xi,\eta)}\nu_-(d\xi)}
    {\int_{X_+}\int_{X_-}e^{-W_k(\xi, \bar\eta)}\nu_-(d\xi)\nu(d\bar\eta)}
    \nu(d\eta).
\end{equation}
Thus, since $\{W_{[N]}\}_{N\in\N}$ (for the definition, see the proof of the first part) is a subsequence of $\{W_k\}_{k\in\N}$,

\begin{equation}\label{f_+ 1_{[1_0^n]}}
    \int_{X_+} f_+ \mathds{1}_{[1_0^n]} d\nu 
    =
    \lim_{N\to\infty}
    \int_{[1_0^n]}
    \frac
    {\int_{X_-}e^{-W_{[N]}(\xi,\eta)}\nu_-(d\xi)}
    {\int_{X_+}\int_{X_-}e^{-W_{[N]}(\xi, \bar\eta)}\nu_-(d\xi)\nu(d\bar\eta)}
    \nu(d\eta).
\end{equation}

Now fix $N\in\N$ and $\eta\in X_+$, consider $W_{[N]}$ as a function of $\xi$.
Clearly, it is a local function, thus by the first part of Theorem \ref{GCB and MCB}, for all $\kappa\in\R$,
\begin{equation}\label{GCB to nu_-}
    \int_{X_-} e^{\kappa[ W_{[N]}(\xi,\eta)-\int_{X_-} W_{[N]}(\xi,\eta)\nu_-(d\xi)   ]}\nu_-(d\xi)
    \leq e^{D \kappa^2 ||\underline \delta (W_{[N]}(\cdot,\eta))||_2^2},
\end{equation}
where $D=4(1-\bar{c}(\Phi))^{-2}$.
By the Cauchy-Schwarz inequality, one can also obtain a lower bound for the integral, in fact, 
\begin{equation}\label{GCB lb to nu_-}
    e^{-D \kappa^2 ||\underline \delta (W_{[N]}(\cdot,\eta))||_2^2}
    \leq
    \int_{X_-} e^{\kappa[ W_{[N]}(\xi,\eta)-\int_{X_-} W_{[N]}(\xi,\eta)\nu_-(d\xi)   ]}\nu_-(d\xi).
\end{equation}
For all $s\in\N$, 
$$
\delta_{-s}(W_{[N]}(\cdot,\eta))=2\beta \Big| \sum_{j=0}^N \frac{\eta_j}{(s+j)^{\alpha}}\Big|\leq 2\beta \sum_{j=0}^N \frac{1}{(s+j)^{\alpha}}.
$$
Hence for $\alpha>\frac{3}{2}$, 
\begin{equation}\label{eq:r-variance}
    ||\underline \delta (W_{[N]}(\cdot,\eta))||_2^2\leq 4\beta^2 \sum_{s=1}^{\infty}\Big(\sum_{j=s}^{\infty}\frac{1}{j^{\alpha}} \Big)^2
    =:
    4\beta^2 C_1(\alpha)<\infty.
\end{equation}
Hence one obtains from (\ref{GCB to nu_-}) and (\ref{GCB lb to nu_-}) that 
\begin{equation}\label{compound ineq for nu_- numerator}
     C_5^{-1}
     \cdot
     e^{-\int_{X_-} W_{[N]}(\xi,\eta)\nu_-(d\xi)}
    \leq
    \int_{X_-} 
    e^{ -W_{[N]}(\xi,\eta)}\nu_-(d\xi)
    \leq
    C_5 
     \cdot
     e^{-\int_{X_-} W_{[N]}(\xi,\eta)\nu_-(d\xi)},
\end{equation}
where $C_5:=e^{4 D \beta^2 C_1(\alpha)}$.

By applying the Gaussian Concentration Bounds to $W_{[N]}$ as a function of both $\xi$ and $\eta$, one can obtain an analogue of (\ref{compound ineq for nu_- numerator}) for $\nu^{(0)}$.
In fact, by using the Cauchy-Schwarz inequality, one can easily check that
\begin{equation}\label{1st compound ineq for nu_0 denominator}
    e^{ -D||\underline\delta(W_{[N]})||^2_2}
    \cdot
    e^{\int_X-W_{[N]}d\nu^{(0)}}
    \leq
    \int_{X}e^{-W_{[N]}}d\nu^{(0)}
    \leq 
    e^{ D||\underline\delta(W_{[N]})||^2_2}
    \cdot
    e^{\int_X-W_{[N]}d\nu^{(0)}}.
\end{equation}
For any $s\in\Z$, 
\begin{equation}
    \delta_s (W_{[N]})
    \leq
    \frac{8\beta}{(|s|+1)^\alpha}+2\beta\sum_{i={|s|+1}}^\infty \frac{1}{i^{\alpha}}
    \leq
    \frac{\beta C_6}{(|s|+1)^{\alpha-1}},
\end{equation}
where $C_6>0$ is only dependent on $\alpha$.
Hence, since $\alpha>3/2$,
\begin{equation}
    ||\underline\delta(W_{[N]})||^2_2
    =
    \sum_{s\in\Z} (\delta_s (W_{[N]}))^2
    \leq
    (\beta C_6)^2\sum_{s\in\Z}
    \frac{1}{(|s|+1)^{2\alpha-2}}=:C_7(\alpha, \beta)<\infty.
\end{equation}
Then one obtains from (\ref{1st compound ineq for nu_0 denominator}) that 
\begin{equation}\label{compound ineq for nu_0 denominator}
    C_8^{-1}
    \cdot
    e^{\int_X-W_{[N]}d\nu^{(0)}}
    \leq
    \int_{X}e^{-W_{[N]}}d\nu^{(0)}
    \leq 
    C_8
    \cdot
    e^{\int_X-W_{[N]}d\nu^{(0)}},
\end{equation}
here $C_8=C_8(\alpha, \beta, h):=D\cdot C_7(\alpha,\beta)$.
By combining (\ref{f_+ 1_{[1_0^n]}}), (\ref{compound ineq for nu_- numerator}) and (\ref{compound ineq for nu_0 denominator}), one concludes that for any $N\in\N$ and $\eta\in X_+$, 
\begin{eqnarray}
    C_9^{-1}
    \cdot
    e^{\int_X W_{[N]}d\nu^{(0)}-\int_{X_-} W_{[N]}(\xi,\eta)\nu_-(d\xi)}
    &\leq&
    \frac{\int_{X_-}e^{-W_{[N]}(\xi,\eta)}\nu_-(d\xi)}
    {\int_{X}e^{-W_{[N]}}d\nu^{(0)}}\\
    &\leq&
    C_9
    \cdot
    e^{\int_X W_{[N]}d\nu^{(0)}-\int_{X_-} W_{[N]}(\xi,\eta)\nu_-(d\xi)},
\end{eqnarray}
here $C_9=C_9(\alpha, \beta, h):=\max\{C_5, C_8\}$.
Hence, by Jensen's inequality,
\begin{multline}\label{apply Jensen's ineq}
    \frac{1}{\nu([1_0^n])}\int_{[1_0^n]}\frac{\int_{X_-}e^{-W_{[N]}(\xi,\eta)}\nu_-(d\xi)}
    {\int_{X}e^{-W_{[N]}}d\nu^{(0)}}\nu(d\eta)
    \geq
    \frac{C_9^{-1}}{\nu([1_0^n])}
    \cdot
    \int_{[1_0^n]}e^{\int_X W_{[N]}d\nu^{(0)}-\int_{X_-} W_{[N]}(\xi,\eta)\nu_-(d\xi)} \nu(d\eta)\\
    \geq
    C_9^{-1}
    \cdot
    \exp{
    \Big(
    \frac{1}{\nu([1_0^n])}
    \int_{[1_0^n]} 
    \Big[
    \int_X  W_{[N]}(\xi,\bar\eta)\nu^{(0)}(d\xi,d\bar\eta)-\int_{X_-} W_{[N]}(\xi,\eta)
    \nu_-(d\xi)
    \Big]\nu(d\eta)
    \Big)}
\end{multline}
Note that for any $\eta\in X_+$, one has
\begin{multline}
\int_X  W_{[N]}(\xi,\bar\eta)\nu^{(0)}(d\xi,d\bar\eta)-\int_{X_-} W_{[N]}(\xi,\eta)
    \nu_-(d\xi)
=
\beta\sum_{i=1}^N\sum_{j=0}^N\frac{\nu_-(\sigma_{-i})[\eta_j-\nu(\sigma_j)]}{(i+j)^\alpha}
\end{multline}
Hence
\begin{multline}\label{lb for 1_0^n integral of diff of means}
\int_{[1_0^n]} 
    \Big[
    \int_X  W_{[N]}(\xi,\bar\eta)\nu^{(0)}(d\xi,d\bar\eta)-\int_{X_-} W_{[N]}(\xi,\eta)
    \nu_-(d\xi)
    \Big]\nu(d\eta)
\\
=
\beta\int_{[1_0^n]} \sum_{i=1}^N\sum_{j=0}^N\frac{\nu_-(\sigma_{-i})[\eta_j-\nu(\sigma_j)]}{(i+j)^\alpha} \nu(d\eta)\\
=
\beta \sum_{i=1}^N\sum_{j=0}^N \frac{\nu_-(\sigma_{-i})}{(i+j)^\alpha} 
\int_{[1_0^n]}[\sigma_j-\nu(\sigma_j)]d\nu.
\end{multline}
Then since both $\sigma_j$ and $\mathds{1}_{[1_0^n]}$ are non-decreasing functions and $\nu$ is a positively correlated measure, by the FKG inequality, for any $j\in\Z_+$, 
\begin{equation}
    \int_{[1_0^n]}[\sigma_j-\nu(\sigma_j)]d\nu
    =
    \int_{X_+}\sigma_j \mathds{1}_{[1_0^n]}-\int_{X_+}\sigma_j d\nu \int_{X_+}\mathds{1}_{[1_0^n]}d\nu\geq 0.
\end{equation}
Furthermore, for $1\leq j\leq n$, one has
\begin{equation}
    \int_{[1_0^n]}[\sigma_j-\nu(\sigma_j)]d\nu
    =
    \nu([1_0^n])(1-\nu(\sigma_j)).
\end{equation}
We also note that by (\ref{eq: GKS-1}), for all $i\in\N$,
\begin{equation}
    \nu_-(\sigma_{-i})\geq 0.
\end{equation}
By combining these arguments, 
\begin{eqnarray}
\beta \sum_{i=1}^N\sum_{j=0}^N \frac{\nu_-(\sigma_{-i})}{(i+j)^\alpha} 
\int_{[1_0^n]}[\sigma_j-\nu(\sigma_j)]d\nu
&\geq&
\beta \sum_{j=0}^n \sum_{i=1}^N \frac{\nu_-(\sigma_{-i})}{(i+j)^\alpha} 
\int_{[1_0^n]}[\sigma_j-\nu(\sigma_j)]d\nu\\
&=&
\beta \sum_{j=0}^n \sum_{i=1}^N \frac{\nu_-(\sigma_{-i})}{(i+j)^\alpha} \nu([1_0^n])(1-\nu(\sigma_j)).
\end{eqnarray}
Hence by (\ref{lb for 1_0^n integral of diff of means}),
\begin{multline}\label{lb for diff of means which diverge as n to infty}
    \frac{1}{\nu([1_0^n])}
    \int_{[1_0^n]} 
    \Big[
    \int_X  W_{[N]}(\xi,\bar\eta)\nu^{(0)}(d\xi,d\bar\eta)-\int_{X_-} W_{[N]}(\xi,\eta)
    \nu_-(d\xi)
    \Big]\nu(d\eta)\\
    \geq
    \frac{\beta}{\nu([1_0^n])}
    \sum_{j=0}^n \sum_{i=1}^N \frac{\nu_-(\sigma_{-i})(1-\nu(\sigma_j))\nu([1_0^n])}{(i+j)^\alpha}\\
    =
    \beta
    \sum_{j=0}^n \sum_{i=1}^N \frac{\nu_-(\sigma_{-i})(1-\nu(\sigma_j))}{(i+j)^\alpha}
\end{multline}
Then (\ref{apply Jensen's ineq}) and (\ref{lb for diff of means which diverge as n to infty}) yields that for any $N\in\N$, 
\begin{multline}
    \frac{1}{\nu([1_0^n])}\int_{[1_0^n]}\frac{\int_{X_-}e^{-W_{[N]}(\xi,\eta)}\nu_-(d\xi)}
    {\int_{X}e^{-W_{[N]}}d\nu^{(0)}}\nu(d\eta)
    \geq
    C_9^{-1}
    \cdot
    \exp{
    \Big(
    \beta
    \sum_{j=0}^n \sum_{i=1}^N \frac{\nu_-(\sigma_{-i})(1-\nu(\sigma_j))}{(i+j)^\alpha}
    \Big)
    }
\end{multline}
Hence (\ref{f_+ 1_{[1_0^n]}}) implies
\begin{eqnarray}
    \frac{1}{\nu([1_0^n])}\int_{[1_0^n]} f_+ d\nu 
    \notag &\geq&
    \lim_{N\to\infty}
    C_9^{-1}
    \cdot
    \exp{
    \Big(
    \beta
    \sum_{j=0}^n \sum_{i=1}^N \frac{\nu_-(\sigma_{-i})(1-\nu(\sigma_j))}{(i+j)^\alpha}
    \Big)
    }\\
    &=&\label{imp lb for the integral av of f_+}
    C_9^{-1}
    \cdot
    \exp{
    \Big(
    \beta
    \sum_{j=0}^n \sum_{i=1}^\infty \frac{\nu_-(\sigma_{-i})(1-\nu(\sigma_j))}{(i+j)^\alpha}
    \Big)
    }.
\end{eqnarray}
It follows from (\ref{eq: speed of convergence of t_n}) that $\lim_{i\to\infty}\nu_-(\sigma_{-i})=\lim_{j\to\infty}\nu(\sigma_j)=\mu(\sigma_0)>0$. 
Thus $\tilde{\kappa}:=\sup_{j\in\Z_+} \nu(\sigma_j)<1$ and there exists $R\in\N$ which depends only on $\nu^{(0)}=\nu_-\times\nu$ and $\mu$ (thus $R$ only depends on the model parameters $\alpha, \beta, h$) such that $\nu_-(\sigma_{-i})\geq \frac{\mu(\sigma_0)}{2}$ for all $i\geq R$.
Hence we obtain from (\ref{imp lb for the integral av of f_+}) that 
\begin{equation}\label{more explicit lb for int av for f_+}
    \frac{1}{\nu([1_0^n])}\int_{[1_0^n]} f_+ d\nu 
    \geq 
    C_9^{-1}
    \cdot
    \exp{
    \Big(
    \frac{\beta (1-\tilde{\kappa})\mu(\sigma_0)}{2}
    \sum_{j=0}^n \sum_{i=R}^\infty \frac{1}{(i+j)^\alpha}
    \Big)
    }.
\end{equation}
One can readily check that the sum in the right-hand side of (\ref{more explicit lb for int av for f_+}) diverges as $n\to\infty$.
In fact, there exists $C_{10}=C_{10}(\alpha)\in (0,1)$ such that $\sum_{i=R}^\infty \frac{1}{(i+j)^\alpha}\geq C_{10} (R+j)^{1-\alpha}$,
therefore, 
$$
\sum_{j=0}^\infty \sum_{i=R}^\infty \frac{1}{(i+j)^\alpha}
\geq
C_{10} \sum_{j=0}^\infty \frac{1}{(R+j)^{\alpha-1}}
=
\infty.
$$
Hence one indeed concludes that 
\begin{equation*}
\sup_{n\in\N}\frac{1}{\nu([1_0^n])}\int_{[1_0^n]} f_+ d\nu 
=
+\infty.
\end{equation*}

\section*{Acknowledgements}
The author expresses gratitude to Aernout van Enter, Evgeny Verbitskiy, and Roberto Fernández for their valuable and insightful discussions. Additionally, the author sincerely appreciates Aernout van Enter's helpful suggestions to improve the manuscript.

\end{document}